\documentclass[12pt]{article}
\usepackage{graphicx}
\usepackage{amsmath,amsthm,amssymb,enumerate}
\usepackage{euscript,mathrsfs}
\usepackage{color}
\usepackage{dsfont}
\usepackage[left=2cm,right=2cm,top=3.5cm,bottom=3.5cm]{geometry}
\usepackage{color}
\usepackage[framemethod=tikz]{mdframed}
\allowdisplaybreaks

\usepackage{soul}

\catcode`\@=11 \@addtoreset{equation}{section}

\catcode`\@=12

\newtheorem{Theorem}{Theorem}[section]
\newtheorem{Proposition}[Theorem]{Proposition}
\newtheorem{Lemma}[Theorem]{Lemma}
\newtheorem{Corollary}[Theorem]{Corollary}

\theoremstyle{definition}
\newtheorem{Definition}[Theorem]{Definition}

\newtheorem{Remark}[Theorem]{Remark}

\newcommand{\bTheorem}[1]{
	\begin{Theorem} \label{T#1} }
	\newcommand{\eT}{\end{Theorem}}

\newcommand{\bProposition}[1]{
	\begin{Proposition} \label{P#1}}
	\newcommand{\eP}{\end{Proposition}}

\newcommand{\bLemma}[1]{
	\begin{Lemma} \label{L#1} }
	\newcommand{\eL}{\end{Lemma}}

\newcommand{\bCorollary}[1]{
	\begin{Corollary} \label{C#1} }
	\newcommand{\eC}{\end{Corollary}}

\newcommand{\bRemark}[1]{
	\begin{Remark} \label{R#1} }
	\newcommand{\eR}{\end{Remark}}

\newcommand{\bDefinition}[1]{
	\begin{Definition} \label{D#1} }
	\newcommand{\eD}{\end{Definition}}

\newcommand{\Del}{\Delta_x}

\newcommand{\Ds}{\mathbb{D}_x}

\newcommand{\vuB}{\vc{u}_B}

\newcommand{\bFormula}[1]{
	\begin{equation} \label{#1}}
	\newcommand{\eF}{\end{equation}}

\newcommand{\Ov}[1]{\overline{#1}}

\newcommand{\aleq}{\stackrel{<}{\sim}}

\newcommand{\vr}{\varrho}

\newcommand{\tvu}{{\tilde \vu}}
\newcommand{\tvt}{\tilde \vt}

\newcommand{\vt}{\vartheta}
\newcommand{\vu}{\vc{u}}

\newcommand{\vc}[1]{{\bf #1}}

\newcommand{\Div}{{\rm div}_x}
\newcommand{\Grad}{\nabla_x}

\newcommand{\dx}{\,{\rm d} {x}}

\newcommand{\dt}{\,{\rm d} t }

\newcommand{\intO}[1]{\int_{\Omega} #1 \ \dx}

\newcommand{\D}{{\rm d}}

\newcommand{\vtB}{\vt_B}

\newcommand{\br}{ \nonumber \\ }

\def\softd{{\leavevmode\setbox1=\hbox{d}%
		\hbox to 1.05\wd1{d\kern-0.4ex{\char039}\hss}}}
\definecolor{Cgrey}{rgb}{0.85,0.85,0.85}
\definecolor{Cblue}{rgb}{0.50,0.85,0.85}
\definecolor{Cred}{rgb}{1,0,0}
\definecolor{fancy}{rgb}{0.10,0.85,0.10}

\newcommand\Cbox[2]{%
	\newbox\contentbox%
	\newbox\bkgdbox%
	\setbox\contentbox\hbox to \hsize{%
		\vtop{
			\kern\columnsep
			\hbox to \hsize{%
				\kern\columnsep%
				\advance\hsize by -2\columnsep%
				\setlength{\textwidth}{\hsize}%
				\vbox{
					\parskip=\baselineskip
					\parindent=0bp
					#2
				}%
				\kern\columnsep%
			}%
			\kern\columnsep%
		}%
	}%
	\setbox\bkgdbox\vbox{
		\color{#1}
		\hrule width  \wd\contentbox %
		height \ht\contentbox %
		depth  \dp\contentbox
		\color{black}
	}%
	\wd\bkgdbox=0bp%
	\vbox{\hbox to \hsize{\box\bkgdbox\box\contentbox}}%
	\vskip\baselineskip%
}

\mdfdefinestyle{MyFrame}{%
	linecolor=black,
	outerlinewidth=1pt,
	roundcorner=5pt,
	innertopmargin=\baselineskip,
	innerbottommargin=\baselineskip,
	innerrightmargin=10pt,
	innerleftmargin=10pt,
	backgroundcolor=white!20!white}

\begin{document}


\title{Conditional regularity for the Navier--Stokes--Fourier system with Dirichlet boundary conditions}

\author{Danica Basari\' c \thanks{The work of D.B., E.F., and H.M. was supported by the
		Czech Sciences Foundation (GA\v CR), Grant Agreement
		21--02411S. The Institute of Mathematics of the Czech Academy of Sciences is supported by RVO:67985840. } 
		\and Eduard Feireisl $^{ *}$
		\and Hana Mizerov\'a $^{ *,\dagger}$
}

\date{}

\maketitle

\medskip

\centerline{$^*$ Institute of Mathematics of the Czech Academy of Sciences}

\centerline{\v Zitn\' a 25, CZ-115 67 Praha 1, Czech Republic}

\bigskip
\centerline{$^\dagger$ Department of Mathematical Analysis and Numerical Mathematics, Comenius University}
\centerline{Mlynsk\' a dolina, 842 48 Bratislava, Slovakia}

\begin{abstract}
	
	We consider the Navier--Stokes--Fourier system with the inhomogeneous boundary conditions for the velocity and the temperature.
	We show that solutions emanating from sufficiently regular data remain regular as long as the 
	density $\vr$, the absolute temperature $\vt$, and the modulus of the fluid velocity $|\vu|$ remain bounded.

\end{abstract}




\section{Introduction}
\label{p}

Standard systems of equations in fluid mechanics including the Navier--Stokes--Fourier system 
governing the motion of a compressible, viscous, and heat conducting fluid are well posed in the 
class of strong solutions on a possibly short time interval $[0,T_{\rm max})$. The recent results of 
Merle at al.  \cite{MeRaRoSz}, \cite{MeRaRoSzbis} strongly indicate that $T_{\rm max}$ may be finite, at least 
in the idealized case of ``isentropic'' viscous flow. Conditional regularity results guarantee that 
a blow up  will not occur as soon as some lower order norms of solutions are controlled. 

We consider the \emph{Navier--Stokes--Fourier system} governing the time evolution of the 
mass density $\vr = \vr(t,x)$, the (absolute) temperature $\vt = \vt(t,x)$, and the velocity 
$\vu = \vu(t,x)$ of a compressible, viscous, and heat conducting fluid:

\begin{mdframed}[style=MyFrame]

\begin{align} 
	\partial_t \vr + \Div (\vr \vu) &= 0, \label{i1} \\
	\partial_t (\vr \vu) + \Div (\vr \vu \otimes \vu) + \Grad p(\vr, \vt) &= \Div \mathbb{S}(\Ds \vu) + \vr \vc{f},\ 
	\Ds \vu = \frac{1}{2} \left( \Grad \vu + \Grad^t \vu \right),  \label{i2} \\
	\partial_t (\vr e(\vr, \vt)) + \Div (\vr e (\vr, 
	\vt) \vu) + \Div \vc{q}(\Grad \vt) &= \mathbb{S} (\Ds \vu) : \Ds \vu - p (\vr, \vt) \Div \vu. \label{i3}
	\end{align}

\end{mdframed}

\noindent
The fluid is Newtonian, the viscous stress $\mathbb{S}$ is given by Newton's rheological law 
\begin{equation} \label{i4}
	\mathbb{S}(\Ds \vu)  = 2\mu \left( \Ds \vu - \frac{1}{3} \Div \vu \mathbb{I} \right) + \eta \Div \vu \mathbb{I},\ \mu > 0,\ \eta \geq 0.
	\end{equation} 
The heat flux obeys Fourier's law 
\begin{equation} \label{i5}
\vc{q}(\Grad \vt) = - \kappa \Grad \vt,\ \kappa > 0.
\end{equation} 
The equation of state for the pressure $p$ and the internal energy $e$ is given by the standard Boyle--Mariotte law of perfect gas, 
\begin{equation} \label{i6}
 p(\vr, \vt) = \vr \vt, \  
	e (\vr, \vt) = c_v \vt,\ c_v > 0.
	\end{equation}
For the sake of simplicity, we suppose that the viscosity coefficients $\mu$, $\eta$, the heat conductivity 
coefficient $\kappa$ as well as the specific heat at constant volume $c_v$ are constant.

There is a large number of recent results concerning conditional regularity for the Navier--Stokes--Fourier system
in terms of various norms. Fan, Jiang, and Ou \cite{FaJiOu} consider a bounded fluid domain $\Omega 
\subset R^3$ with the conservative 
boundary conditions
\begin{equation} \label{i7}
\vu|_{\partial \Omega} = 0, \ \Grad \vt \cdot \vc{n}|_{\partial \Omega} = 0.
\end{equation}
The same problem is studied by Sun, Wang, and Zhang \cite{SuWaZh} and later 
by Huang, Li, Wang \cite{Huang}.  There are results for the Cauchy problem 
$\Omega = R^3$ by Huang and Li \cite{HuaLi}, and Jiu, Wang and Ye \cite{JiuWanYe}. 
Possibly the best result so far has been established in \cite{FeWeZh}, where the blow 
up criterion for both the Cauchy problem and the boundary value problem \eqref{i7} is formulated in terms 
of the maximum of the density and a Serrin type regularity for the temperature: 
\[
\limsup_{t \to T_{\rm max}-} \left( \| \vr(t, \cdot) \|_{L^\infty} + \| \vt - \vt_\infty \|_{L^s(0,t)(L^r)}  \right) = \infty, \ \frac{3}{2} < r \leq \infty,\ 1 \leq s \leq \infty,\ \frac{2}{s} + \frac{3}{r} \leq 2, 
\]
where $\vt_\infty$ denotes the far field temperature in the Cauchy problem, cf. also the previous results by 
Wen and Zhu \cite{WenZhu1}, \cite{WenZhu2}.

Much less is known in the case of the Dirichlet boundary conditions 
\begin{equation} \label{i8}
\vu|_{\partial \Omega} = \vuB,\ 	\vt|_{\partial \Omega} = \vtB.
	\end{equation}
Fan, Zhi, and Zhang \cite{FaZiZh} showed that a strong solution of the Navier--Stokes--Fourier system remains regular up to 
a time $T > 0$ if (i) $\Omega \subset R^2$ is a bounded domain, (ii) $\vuB = 0$, $\vtB = 0$, and (iii) 
\begin{equation} \label{i9}
\limsup_{t \to T-} \left( \| \vr \|_{L^\infty} + \| \vt \|_{L^\infty} \right) < \infty.	
	\end{equation}

All results mentioned above describe fluids in a conservative regime, meaning solutions are close to equilibrium 
in the long run. However, many real world applications concern fluids out of equilibrium driven by possibly 
large driving forces $\vc{f}$ and/or inhomogeneous boundary conditions. The iconic examples 
are the Rayleigh--B\' enard and Taylor--Couette flows where the fluid is driven to a turbulent regime 
by a large temperature gradient and large boundary velocity, respectively, see Davidson \cite{DAVI}. 

Motivated 
by these physically relevant examples, we consider a fluid confined to a bounded domain $\Omega \subset R^3$ 
with \emph{impermeable boundary}, where the temperature and the (tangential) velocity are given on $\partial \Omega$,  
\begin{mdframed}[style=MyFrame]
	\begin{align} 
		\vt|_{\partial \Omega} &= \vtB,\ \vtB = \vtB(x),\ \vtB > 0 \ \mbox{on}\ \partial \Omega, \label{i10} \\
		\vu|_{\partial\Omega} &= \vuB,\ \vuB = \vuB(x),\ \vuB \cdot \vc{n} = 0 
		\ \mbox{on}\ \partial \Omega. \label{i11}	
	\end{align}	
	
\end{mdframed}
\noindent
The initial state of the fluid is prescribed:
\begin{mdframed}[style=MyFrame]
\begin{equation} \label{i12}
\vr(0, \cdot) = \vr_0,\ \vr_0 > 0 \ \mbox{in}\ \Ov{\Omega},\ 
\vt(0, \cdot) = \vt_0,\ \vt_0 > 0 \ \mbox{in}\ \Ov{\Omega},\ \vu(0, \cdot) = \vu_0. 	
	\end{equation}	

	\end{mdframed}
\noindent
The initial and boundary data are supposed to satisfy suitable \emph{compatibility conditions} specified below.

The existence of local in time strong solutions for the problem \eqref{i1}--\eqref{i6}, endowed with the inhomogeneous boundary conditions \eqref{i10}, \eqref{i11} was established by Valli \cite{Vall2}, \cite{Vall1} , see also Valli and Zajaczkowski \cite{VAZA}. The solution exists on a maximal time interval $[0, T_{\rm max})$, 
$T_{\rm max} > 0$. 
Our goal is to show that if $T_{\rm max} < \infty$, then necessarily 
\begin{equation} \label{i13}
\limsup_{t \to T_{\rm max}-} \Big( \| \vr (t, \cdot)  \|_{L^\infty(\Omega)} + 
\| \vt (t, \cdot) \|_{L^\infty(\Omega)} + \| \vu (t, \cdot) \|_{L^\infty(\Omega; R^3)} \Big) = \infty.
\end{equation}
The proof is based on deriving suitable {\it a priori} bounds assuming boundedness of all norms involved 
in \eqref{i13} as well as the norm of the initial/boundary data in a suitable function space. Although approach shares 
some similarity with Fang, Zi, and Zhang \cite{FaZiZh}, essential modifications must be made to accommodate 
the inhomogeneous boundary data as well as the driving force $\vc{f}$. The importance of conditional regularity 
results in numerical analysis of flows with uncertain initial data was discussed recently in \cite{FeiLuk2021}.

The paper is organized as follows. In Section \ref{M}, we introduce the class of strong solutions to the Navier--Stokes--Fourier system and state our main result concerning conditional regularity. The remaining 
part of the paper is devoted to the proof of the main result -- deriving suitable {\it a priori} bounds.
In Section \ref{e} we recall the standard energy estimates that hold even in the class of weak solutions.
Section \ref{g} is the heart of the paper. We establish the necessary estimates on the velocity gradient 
by means of the celebrated Gagliardo--Nirenberg interpolation inequality. In Section \ref{s}, higher order estimates 
on the velocity gradient are derived, and, finally, the estimates are closed by proving bounds on the 
temperature time derivative in Section \ref{d}. This last part borrows the main ideas from \cite{FeNoSun1}.

\section{Strong solutions, main result}
\label{M}

We start the analysis by recalling the concept of strong solution introduced by Valli \cite{Vall1}. 
Similarly to the boundary data $\vuB$, $\vtB$ we suppose that the driving force $\vc{f} = \vc{f}(x)$ is independent 
of time, meaning we deal with an autonomous problem. Following \cite{Vall1}, we suppose that $\Omega \subset R^3$ is a bounded domain with $\partial \Omega$ of class $C^4$. 

We assume the data belong to the following class:
\begin{align}
\vr_0 &\in W^{3,2}(\Omega),\ 0 < \underline{\vr}_0 \leq \min_{x \in \Omega} \vr_0 (x), \br	
\vt_0 &\in W^{3,2}(\Omega),\ 0 <  \underline{\vt}_0\leq \min_{x \in \Omega} \vt_0 (x), \br
\vu_0 &\in W^{3,2}(\Omega; R^3), \br
\vtB &\in W^{\frac{7}{2}} (\partial \Omega),\ 0 <  \underline{\vt}_B \leq \min_{x \in \partial \Omega} \vtB (x), \br
\vuB &\in W^{\frac{7}{2}} (\partial \Omega; R^3),\ \vuB \cdot \vc{n} = 0, \br 
\vc{f}& \in W^{2,2}(\Omega; R^3).	
	\label{M1}
	\end{align}
In addition, the data must satisfy the compatibility conditions 
\begin{align}
\vt_0 = \vtB,\ \vu_0 &= \vuB \ \mbox{on}\ \partial \Omega, \br	 
\vr_0 \vu_0 \cdot \Grad \vu_0 + \Grad p(\vr_0, \vt_0) &= \Div \mathbb{S} (\Ds \vu_0) + \vr_0 \vc{f} \ \mbox{on}\ 
\partial \Omega, \br
\vr_0 \vu_0 \cdot \Grad \vt_0 + \Div \vc{q}(\vt_0) &= \mathbb{S}(\Ds \vu_0) : \Ds \vu_0 - p(\vr_0, \vt_0) \Div 
\vu_0 \ \mbox{on}\ \partial \Omega.
\label{M2}	
	\end{align}
We set 
\begin{equation} \label{M3}
\mathcal{D}_0 = \max\left\{ \| (\vr_0, \vt_0, \vu_0) \|_{W^{3,2}(\Omega; R^5)}, \frac{1}{\underline{\vr}_0},\ 
\frac{1}{\underline{\vt}_0}, \frac{1}{\underline{\vt}_B},  \| \vtB \|_{W^{\frac{7}{2}}(\partial \Omega)}, \| \vuB \|_{W^{\frac{7}{2}}(\partial \Omega; R^3)} , \| \vc{f} \|_{W^{2,2}(\Omega; R^3)}	  \right\}.
\end{equation}

\subsection{Local existence}

The following result was proved by Valli \cite[Theorem A]{Vall1} (see also \cite{Vall2}). 

\begin{Theorem} \label{TVal} {\bf (Local existence of strong solutions)}
Let $\Omega \subset R^3$ be a bounded domain of class $C^4$. Suppose that the data $(\vr_0, \vt_0, \vu_0)$, 
$(\vtB, \vuB)$ and $\vc{f}$ belong to the class  \eqref{M1} and satisfy the compatibility conditions \eqref{M2}. 

Then there exists a maximal time $T_{\rm max} > 0$ such that the Navier--Stokes--Fourier system \eqref{i1}--\eqref{i6}, 
with the boundary conditions \eqref{i10}, \eqref{i11}, and the initial conditions \eqref{i12} admits a 
solution $(\vr, \vt, \vu)$ in $[0, T_{\rm max}) \times \Omega$ unique in the class
\begin{align} 
\vr,\ \vt &\in C([0,T]; W^{3,2}(\Omega)),\ \vu \in C([0,T]; W^{3,2}(\Omega; R^3)), \br
\vt &\in L^2(0,T; W^{4,2}(\Omega)),\ \vu \in L^2(0,T; W^{4,2}(\Omega; R^3)) 
		\label{M4}
	\end{align}
for any $0 < T < T_{\rm max}$. The existence time $T_{\rm max}$ is bounded below by a quantity $c(\mathcal{D}_0)$ 
depending solely on the norms of the data specified in \eqref{M3}. In particular, 
\begin{equation} \label{M5}
\lim_{\tau \to T_{\rm max}-} \| (\vr, \vt, \vu) (\tau, \cdot) \|_{W^{3,2}(\Omega; R^5)} = \infty.	
	\end{equation}   	
	\end{Theorem}

\subsection{Blow up criterion, conditional regularity}
\label{c}

Our goal is to show the following result. 

\begin{mdframed}[style=MyFrame]
	
	\begin{Theorem} \label{MT} {\bf (Blow up criterion)}
	Under the hypotheses of Theorem \ref{TVal}, suppose that the maximal existence time $T_{\rm max} < \infty$ is finite. 
	
	Then 
	\begin{equation} \label{M6}
	\limsup_{\tau \to T_{\rm max}-} \left\| (\vr, \vt, \vu)(\tau, \cdot) \right\|_{L^\infty (\Omega; R^5)} = \infty.	
		\end{equation}	
		
		\end{Theorem}
	
	\end{mdframed}

Theorem \ref{MT} is in the spirit of the blow up criteria for general parabolic systems -- the solution remains 
regular as long as it is bounded. Of course, our problem in question is of mixed hyperbolic--parabolic type.

The proof of Theorem \ref{MT} follows from suitable {\it a priori} bounds applied on a compact time interval. 

\begin{Proposition} \label{PT} {\bf (Conditional regularity)}
	
\noindent	
Under the hypotheses of Theorem \ref{TVal}, let $(\vr, \vt, \vu)$ be the strong solution of the Navier--Stokes--Fourier system belonging to the class \eqref{M4} and satisfying	
\begin{equation} \label{c1}
	\sup_{(\tau,x) \in [0,T) \times \Omega} \vr (\tau,x) \leq \Ov{\vr},\ \sup_{(\tau,x) \in [0,T) \times \Omega} \vt (\tau,x) \leq \Ov{\vt},\ 
	\sup_{(\tau ,x) \in [0,T) \times \Omega} | \vu (\tau,x) | \leq \Ov{u}
\end{equation}	
for some $T < T_{\rm max}$. 

Then there is a quantity $c(T, \mathcal{D}_0, \Ov{\vr}, \Ov{\vt}, \Ov{u})$, bounded for bounded arguments, such that 
\begin{equation} \label{M7}
\sup_{\tau \in [0,T)} \max \left\{ \| (\vr, \vt, \vu) (\tau, \cdot) \|_{W^{3,2}(\Omega; R^5)} ; 
\sup_{x \in \Omega} \frac{1}{\vr (\tau,x) } ; \sup_{x \in \Omega} \frac{1}{\vt (\tau,x) } \right\}	
\leq c(T, \mathcal{D}_0, \Ov{\vr}, \Ov{\vt}, \Ov{u}).
	\end{equation}
	
	\end{Proposition}
	
In view of Theorem \ref{TVal}, the conclusion of Theorem \ref{MT} follows from Proposition \ref{PT}. 
The rest of the paper is therefore devoted to the proof of Proposition \ref{PT}.

\begin{Remark} \label{Rem1}
	
	As observed in \cite{FeiLuk2022}, the conditional regularity results established in Proposition \ref{PT} 
	gives rise to \emph{stability} with respect to the data. More specifically, the maximal existence time 
	$T_{\rm max}$ is a lower semicontinuous function of the data with respect to the topologies in \eqref{M1}.

	\end{Remark}

\begin{Remark} \label{Rem2}
	
	Conditional regularity results in combination with the weak--strong uniqueness principle in the class of measure--valued solutions is an efficient tool for proving convergence of numerical schemes, see 
	\cite[Chapter 11]{FeLMMiSh}. The concept of measure--valued solutions to the Navier--Stokes--Fourier system with inhomogeneous Dirichlet boundary conditions has been introduced recently by Chaudhuri \cite{Chaudh}.

	\end{Remark}	

\section{Energy estimates}
\label{e}

To begin, it is suitable to extend the boundary data into $\Omega$. For definiteness, we consider the (unique) solutions 
of the Dirichlet problem 
\begin{equation} \label{e5}
	\begin{aligned}
	\Del \tvt &= 0 \ \mbox{in}\ \Omega,\ \tvt|_{\partial \Omega} = \vtB, \\
	\Div \mathbb{S}(\Ds \tvu) &= 0 \ \mbox{in}\ \Omega,\ \tvu|_{\partial \Omega} = \vuB.
	\end{aligned}
\end{equation}
By abuse of notation, we use the same symbol $\vtB$, $\vuB$ for both the boundary values and their $C^1$ extensions $\tvt = \tvt(x)$, $\tvu = \tvu(x)$ inside $\Omega$.

We start with the ballistic energy equality, see \cite[Section 2.4]{ChauFei},
\begin{align} 
	\frac{\D }{\dt} &\intO{ \left( \frac{1}{2} \vr |\vu - \vuB|^2 + \vr e - \vtB \vr s \right) }
	+\intO{ \frac{\vtB}{\vt} \left( \mathbb{S}(\Ds \vu) : \Ds \vu + \kappa \frac{ |\Grad \vt|^2 }{\vt} \right) } 	\br
	&= - \intO{ \Big(  \vr \vu \otimes \vu  + p \mathbb{I} - \mathbb{S} (\Ds \vu)  \Big) : \Ds \vuB } 
	+ \frac{1}{2} \intO{ \vr \vu \cdot \Grad |\vuB|^2 } \br 
	&+ \intO{ \vr (\vu - \vuB) \cdot \vc{f} } - \intO{ \vr s \vu \cdot \Grad \vtB } + \kappa \intO{ \frac{\Grad \vt}{\vt} \cdot \Grad \vtB },  \label{e1}
\end{align}
where we have introduced the entropy
\[
s = c_v \log (\vt) - \log(\vr) .
\]
Thus the choice \eqref{e5} yields the following bounds
\begin{align}
	\sup_{t \in [0,T) } \intO{ \vr | \log(\vt) | (t, \cdot) } \leq c(T, \mathcal{D}_0, \Ov{\vr}, \Ov{\vt}, \Ov{u} ), \label{e2}\\
	\int_0^T \intO{ |\Grad \vu |^2 } \dt  \leq C(\Ov{\vr}, \Ov{\vt}, \Ov{u}; {\rm data} ) 
	\ \Rightarrow \ \int_0^T \| \vu \|^2_{W^{1,2}(\Omega; R^3) } \dt \leq c(T, \mathcal{D}_0, \Ov{\vr}, \Ov{\vt}, \Ov{u} ), \label{e3}\\
	\int_0^T \intO{ \left( |\Grad \vt |^2 + |\Grad \log(\vt) |^2 \right) } \dt  \leq c(T, \mathcal{D}_0, \Ov{\vr}, \Ov{\vt}, \Ov{u} ), \br \Rightarrow \ \int_0^T \| \vt \|^2_{W^{1,2}(\Omega) } \dt + 
	\int_0^T \| \log (\vt) \|^2_{W^{1,2}(\Omega) } \dt \leq c(T, \mathcal{D}_0, \Ov{\vr}, \Ov{\vt}, \Ov{u} ) .\label{e4}
	\end{align}  
	
\section{Estimates of the velocity gradient}
\label{g}

This section is the heart of the paper. In principle, we follow the arguments similar to Fang, Zi, and Zhang \cite[Section 3]{FaZiZh}
but here adapted to the inhomogeneous boundary conditions.

\subsection{Estimates of the velocity material derivative}

Let us introduce the material derivative of a function $g$, 
\[
D_t g = \partial_t g + \vu \cdot \Grad g.
\]
Accordingly, we may rewrite the momentum equation \eqref{i2} as 
\begin{equation} \label{g1}
	\vr D_t \vu + \Grad p = \Div \mathbb{S} + \vr \vc{f}.	
	\end{equation}
Now, consider the scalar product 
of the momentum equation \eqref{g1} with $D_t (\vu - \vuB)$, 
\begin{equation} \label{g1b}
\vr |D_t \vu|^2 + \Grad p \cdot D_t (\vu - \vuB) = \Div \mathbb{S}(\Ds \vu) \cdot D_t (\vu - \vu_B) + \vr \vc{f} \cdot D_t (\vu- \vuB) + \vr D_t \vu \cdot D_t \vuB.
\end{equation}

The next step is integrating \eqref{g1b} over $\Omega$. Here and hereafter we use the hypothesis
$\vuB \cdot \vc{n}|_{\partial \Omega} = 0$ yielding
\begin{equation} \label{g1a}
	D_t (\vu - \vu_B) |_{\partial \Omega} = \left( \partial_t \vu - \vu \cdot \Grad (\vu - \vu_B) \right) 
	|_{\partial \Omega} = - \vuB \cdot \Grad (\vu - \vuB)|_{\partial \Omega} = 0.
	\end{equation}
Writing 
\[
\Div \mathbb{S}(\Ds \vu) = \mu \Del \vu + \left( \eta+ \frac{\mu}{3} \right) \Grad \Div \vu,
\]
and making use of \eqref{g1a} we obtain
\begin{align}
&\intO{ \Div \mathbb{S}(\Ds \vu) \cdot D_t (\vu - \vu_B) }   \br
=& - \intO{ \mathbb{S}(\Ds \vu)  : \Grad \partial_t \vu } \br
&- \mu \intO{ \Grad \vu : \Grad \big(\vu \cdot \Grad (\vu - \vuB) \big) } - 
\left( \eta+ \frac{\mu}{3} \right) \intO{ \Div \vu \  \Div \big(\vu \cdot \Grad (\vu - \vuB)\big) } \br
=& - \frac{1}{2} \frac{\D }{\dt} \intO{ \mathbb{S} (\Ds \vu) : \Ds \vu } \br 
&-\mu \intO{ \Grad \vu : \Grad \big(\vu \cdot \Grad (\vu - \vuB) \big) } - 
\left( \eta+\frac{\mu}{3} \right) \intO{ \Div \vu \ \Div \big(\vu \cdot \Grad (\vu - \vuB) \big) },
\label{g1g}
\end{align}
where, furthermore, 
\begin{align}
\intO{\Grad \vu : \Grad (\vu \cdot \Grad \vu)  }  &= 
\intO{\Grad \vu : (  \Grad \vu \cdot \Grad \vu)  } + \frac{1}{2} \intO{ \vu \cdot \Grad |\Grad \vu|^2 }  \br 
&= \intO{\Grad \vu : (  \Grad \vu \cdot \Grad \vu)  } - \frac{1}{2} \intO{ \Div \vu |\Grad  \vu|^2 }
\label{g1c}
\end{align}
Note carefully we have used $\vu \cdot \vc{n}|_{\partial \Omega} = 0$ in the last integration. Similarly, 
\begin{align}
\intO{ \Div \vu \ \Div (\vu \cdot \Grad \vu) }	= 
\intO{ \Div \vu \ \Grad \vu : \Grad^t \vu } - \frac{1}{2} \intO{ (\Div \vu )^3 }.
\label{g1ca}
	\end{align}
Thus summing up the previous observations, we get 
\begin{align}
\frac{1}{2} \frac{\D }{\dt} \intO{ \mathbb{S} (\Ds \vu) : \Ds \vu } &+ 
\frac{1}{2} \intO{ \vr |D_t \vu|^2 } + \intO{ \Grad p \cdot D_t (\vu - \vuB) } \br &\leq c(T, \mathcal{D}_0, \Ov{\vr}, \Ov{\vt}, \Ov{u} ) \left(1  + \intO{ |\Grad \vu |^3 } \right).	
\label{g2}
	\end{align}

Moreover, 
\begin{align}
\intO{ \Grad p \cdot D_t (\vu - \vuB) }  &= - \intO{ p \ \Div ( D_t (\vu - \vuB) ) }	\br
&= - \intO{ p \  \Div  D_t \vu  } + \intO{ p \ \Div (\vu \cdot \Grad \vuB)}, \label{g3a}
	\end{align}
where 
\begin{align*}
p \ \Div D_t \vu  &= \partial_t (p \  \Div \vu) - \big( \partial_t p + \Div (p \vu) \big) \Div \vu 
+ \Div (p \vu) \Div \vu + p \ \Div (\vu \cdot \Grad \vu) \\[0.1cm] 
&= \partial_t (p \ \Div \vu) - \big( \partial_t p + \Div (p \vu) \big) \Div \vu + 
p \Grad \vu : \Grad^t \vu + \Div \big( p\vu \ \Div \vu \big).
\nonumber
\end{align*}
As $\vu \cdot \vc{n}|_{\partial \Omega} = 0$, we have 
\[
\intO{ \Div \big( p \vu \ \Div \vu \big) } = 0,
\]
and the above estimates together with \eqref{g2} give rise to 
\begin{align}
	\frac{1}{2} \frac{\D }{\dt} \intO{ \mathbb{S} (\Ds \vu) : \Ds \vu } &- \frac{\D }{\dt} \intO{ p \Div \vu } + 
	\frac{1}{2} \intO{ \vr |D_t \vu|^2 }  \br
	&\leq c(T, \mathcal{D}_0, \Ov{\vr}, \Ov{\vt}, \Ov{u} ) \left(1  + \intO{ |\Grad \vu |^3 } \right) -\intO{ \big( \partial_t p + \Div (p \vu) \big) \Div \vu }.	
	\nonumber
\end{align}

Finally, we realize 
\[
\partial_t p + \Div (p \vu) = \vr D_t \vt
\]
to conclude 
\begin{align}
	\frac{1}{2} \frac{\D }{\dt} \intO{ \mathbb{S} (\Ds \vu) : \Ds \vu } &- \frac{\D }{\dt} \intO{ p \Div \vu } + 
	\frac{1}{2} \intO{ \vr |D_t \vu|^2 } \br 
		&\leq c(T, \mathcal{D}_0, \Ov{\vr}, \Ov{\vt}, \Ov{u} ) \left(1  + \intO{ \vr |D_t \vt| |\Grad \vu|    } + \intO{ |\Grad \vu |^3 } \right).	
	\label{g2a}
\end{align}

\subsection{Higher order velocity material derivative estimates}

Following \cite[Section 3, Lemma 3.3]{FaZiZh}, see also Hoff \cite{HOF1}, we deduce 
\begin{align} 
\vr &D^2_t \vu + \Grad \partial_t p + \Div (\Grad p \otimes \vu)	 \br &= \mu \Big( \Del \partial_t \vu + \Div (\Del \vu \otimes \vu) \Big) + \left( \eta + \frac{\mu}{3} \right) \Big(
\Grad \Div \partial_t \vu + \Div \left( (\Grad \Div \vu) \otimes \vu \right) \Big) + \vr \vu \cdot \Grad \vc{f}.
\label{g3}
	\end{align}
Next, we compute
\begin{align}
D_t \vuB = \vu \cdot \Grad \vuB, \quad 
D^2_t \vuB &= \partial_t \vu \cdot \Grad \vuB + \vu \cdot \Grad (\vu \cdot \Grad \vuB) \br &= 
D_t \vu \cdot \Grad \vuB - ( \vu \cdot \Grad \vu ) \cdot \Grad \vuB +  \vu \cdot \Grad (\vu \cdot \Grad \vuB) \br
&= D_t \vu \cdot \Grad \vuB + (\vu \otimes \vu) : \nabla^2_x \vuB.
\label{g3b}
\end{align} 
Consequently, we may rewrite \eqref{g3} in the form 
\begin{align} 
	\vr &D^2_t ( \vu - \vuB)  + \Grad \partial_t p + \Div (\Grad p \otimes \vu)	 \br &= \mu \Big( \Del \partial_t \vu + \Div (\Del \vu \otimes \vu) \Big) + \left( \eta + \frac{\mu}{3} \right) \Big(
	\Grad \Div \partial_t \vu + \Div \left( (\Grad \Div \vu) \otimes \vu \right) \Big) + \vr \vu \cdot \Grad \vc{f}
	\br &- \vr D_t \vu \cdot \Grad \vuB - \vr (\vu \otimes \vu) : \nabla^2_x \vuB.
	\label{g3c}
\end{align}

The next step is considering the scalar product of \eqref{g3c} with $D_t (\vu - \vuB)$
and integrating over $\Omega$. The resulting 
integrals can be handled as follows:
\begin{align}
\vr D^2_t (\vu - \vuB) \cdot D_t (\vu - \vuB) &=
\vr \frac{1}{2} D_t | D_t (\vu - \vuB) |^2 \br
&= \frac{1}{2} \vr \left( \partial_t | D_t (\vu - \vuB) |^2 + \vu \cdot \Grad | D_t (\vu - \vuB) |^2 \right)\br
&= \frac{1}{2} \partial_t \left( \vr | D_t (\vu - \vuB) |^2 \right) + 
\frac{1}{2} \Div \left( \vr \vu | D_t (\vu - \vuB) |^2 \right),
\nonumber
\end{align}
where we have used the equation of continuity \eqref{i1}.
Seeing that $\vu \cdot \vc{n}|_{\partial \Omega} = 0$ we get 
\begin{equation} \label{G1}
\intO{ \vr D^2_t (\vu - \vuB) \cdot D_t (\vu - \vuB) } = 
\frac{\D }{\dt} \frac{1}{2} \intO{ \vr |D_t (\vu - \vuB) |^2 }. 
\end{equation}

Similarly, 
\begin{align}
&\intO{ \Big( \Grad \partial_t p + \Div (\Grad p \otimes \vu) \Big) \cdot D_t (\vu - \vuB) } \br 
&=- 	\intO{ \Big( \partial_t p + \Div (p \vu) \Big) \Div D_t (\vu - \vuB) } \br 
&+ \intO{ \Big( \Div (p\vu) \Div D_t (\vu - \vuB) - \Grad p \otimes \vu : \Grad D_t (\vu - \vuB)  \Big) }, 
\label{G2}	
	\end{align}	
where 
\begin{align}
&\intO{ \Grad p \otimes \vu : \Grad D_t (\vu - \vuB) } \br
&= - \intO{ p \Grad \vu : \Grad D_t (\vu - \vuB) } + 
\intO{ \Grad ( p \vu ) : \Grad D_t (\vu - \vuB)  }.	
	\nonumber
	\end{align}
In addition, as $D_t (\vu - \vuB)$ vanishes on $\partial \Omega$, we can perform by parts integration in the 
last integral obtaining 
\[
\intO{ \Grad ( p \vu ) : \Grad D_t (\vu - \vuB)  } = \intO{ \Div ( p \vu )  \Div D_t (\vu - \vuB)  }.
\]
Thus, similarly to the preceding section, we conclude 
\begin{align}
	&\intO{ \Big( \Grad \partial_t p + \Div (\Grad p \otimes \vu) \Big) \cdot D_t (\vu - \vuB) } \br 
	&= - 	\intO{ \vr D_t \vt \Div D_t (\vu - \vuB) } \
	+ \intO{ p \Grad \vu : \Grad D_t (\vu - \vuB) }.
	\label{G3}	
\end{align}	
Analogously, 
\begin{align}
&\intO{ \Big( \Del \partial_t \vu + \Div (\Del \vu \otimes \vu) \Big) \cdot D_t (\vu - \vuB) }\br	
&= - \intO{ \Grad \partial_t \vu : \Grad  D_t (\vu - \vuB) }
- \intO{ ( \Del \vu \otimes \vu ) : \Grad D_t (\vu - \vuB) } \br 
&= - \intO{ \Grad D_t \vu : \Grad  D_t (\vu - \vuB)  }  - \intO{ \Big(  \Del \vu \otimes \vu - 
	\Grad ( \vu \cdot \Grad \vu) \Big) : \Grad D_t (\vu - \vuB) }, 
	\label{G4}
	\end{align}
where, using summation convention,  
\begin{align}
&\intO{  \big(\Del \vu \otimes \vu\big) : \Grad D_t (\vu - \vuB) } \br 
&= \intO{ \partial_{x_k} \Big( u_j\partial_{x_k} u_i  \Big) \partial_{x_j} D_t( \vu - \vuB)_i } 
- \intO{ \partial_{x_k} u_i \partial_{x_k} u_j  \partial_{x_j} D_t( \vu - \vuB)_i } \br 
&= \intO{ \partial_{x_j} \Big( u_j \partial_{x_k} u_i  \Big) \partial_{x_k} D_t( \vu - \vuB)_i } 
- \intO{ \partial_{x_k} u_i \partial_{x_k} u_j  \partial_{x_j} D_t( \vu - \vuB)_i } \br 
&= \intO{ \Div \vu \ \Grad \vu : \Grad  D_t( \vu - \vuB) } \br 
&+
\intO{ \Big( u_j \partial_{x_k} \partial_{x_j} u_i  \Big) \partial_{x_k} D_t( \vu - \vuB)_i } 
- \intO{ \partial_{x_k} u_i \partial_{x_k} u_j  \partial_{x_j} D_t( \vu - \vuB)_i } \br 
&= \intO{ \Grad (\vu \cdot \Grad \vu ) :  \Grad D_t( \vu - \vuB) } + \intO{ \Div \vu \ \Grad \vu : \Grad  D_t( \vu - \vuB) } \br 
 &-\intO{  \partial_{x_j} u_i \partial_{x_k} u_j \partial_{x_k} D_t( \vu - \vuB)_i } 
- \intO{ \partial_{x_k} u_i \partial_{x_k} u_j  \partial_{x_j} D_t( \vu - \vuB)_i }.
\label{G5} 
\end{align}
Summing up \eqref{G4}, \eqref{G5} we conclude 
\begin{align}
	&\intO{ \Big( \Del \partial_t \vu + \Div (\Del \vu \otimes \vu) \Big) \cdot D_t (\vu - \vuB) }\br	
	&= - \intO{ \Grad D_t \vu : \Grad  D_t (\vu - \vuB)  }  -\intO{ \Div \vu \ \Grad \vu : \Grad  D_t( \vu - \vuB) } \br 
	&+\intO{  \partial_{x_j} u_i \partial_{x_k} u_j  \partial_{x_k} D_t( \vu - \vuB)_i } 
	+ \intO{ \partial_{x_k} u_i \partial_{x_k} u_j  \partial_{x_j} D_t( \vu - \vuB)_i }.	
	\label{G6}
\end{align}

Estimating the remaining integrals in \eqref{g3c} in a similar manner we may infer
\begin{align}
\frac{1}{2} \frac{\D}{\dt} &\intO{ \vr |D_t (\vu - \vuB) |^2 }	+ \mu \intO{ |\Grad D_t (\vu - \vuB) |^2 } + \left( \eta +\frac{\mu}{3}\right) \intO{ |\Div D_t ( \vu - \vuB )|^2 } \br 
&\leq  c(T, \mathcal{D}_0, \Ov{\vr}, \Ov{\vt}, \Ov{u} ) \left( 1 + \intO{ \vr |D_t \vt |^2 } + \intO{ |\Grad \vu |^4 }
+ \intO{ \vr |D_t \vu |^2 } \right).
\label{g4}
	\end{align}
cf. \cite[Section 3, Lemma 3.3]{FaZiZh}.

\subsection{Velocity decomposition}

Following the original idea of Sun, Wang, and Zhang \cite{SuWaZh1}, we decompose the velocity field in the form:
\begin{align}
	\vu &= \vc{v} + \vc{w}, \label{g5} \\  
\Div \mathbb{S}(\Ds \vc{v} ) &= \Grad p \ \mbox{in}\ (0,T) \times \Omega ,\ \vc{v}|_{\partial \Omega} = 0, \label{g6} \\
\Div \mathbb{S}(\Ds \vc{w} ) &= \vr D_t \vu - \vr \vc{f} \ \mbox{in}\ (0,T) \times \Omega,\ 
\vc{w}|_{\partial \Omega} = \vuB. \label{g7} 	
	\end{align}

Since 
\[
\Div \mathbb{S}(\Ds \partial_t \vc{v} ) = \Grad \partial_t p \ \mbox{in}\ (0,T) \times \Omega ,\ \vc{v}|_{\partial \Omega} = 0, 
\]
we get 
\begin{equation} \label{g8}
\intO{ \partial_t p \ \Div  \vc{v} } = - \intO{ \Grad \partial_t p \cdot \vc{v} } = \frac{1}{2} \frac{\D }{\dt} \intO{ \mathbb{S}(\Ds \vc{v}) : \Ds \vc{v} }.  
\end{equation} 	

Moreover, the standard elliptic estimates for the Lam\' e operator yield:
\begin{align} 
	\| \vc{v} \|_{W^{1,q}(\Omega; R^3)} &\leq c(q, \Ov{\vr}, \Ov{\vt}) \ \mbox{for all}\ 1 \leq q < \infty, \label{g9} \\ 
	\| \vc{v} \|_{W^{2,q}(\Omega; R^3)} &\leq c(q, \Ov{\vr}, \Ov{\vt}) \left( \| \Grad \vr \|_{L^q(\Omega; R^3)} + \| \Grad \vt \|_{L^q(\Omega; R^3)} \right),\ 
	1 < q < \infty 
	. \label{g10} 
\end{align}
Similarly, 
\begin{equation} \label{g11}
	\| \vc{w} \|_{W^{2,2}(\Omega; R^3)} \leq c(T, \mathcal{D}_0, \Ov{\vr}, \Ov{\vt}, \Ov{u} ) \left( 1 + \| \sqrt{\vr} \partial_t \vu \|_{L^2(\Omega; R^3)} + 
	\| \Grad \vu \|_{L^2(\Omega; R^{3 \times 3})} \right).   
\end{equation}
The estimates \eqref{g9}--\eqref{g11} are uniform in the time interval $[0,T)$.

\subsection{Temperature estimates}

Similarly to Fang, Zi, Zhang \cite[Section 3, Lemma 3.4]{FaZiZh} we multiply the internal energy equation \eqref{i3} on $\partial_t \vt$ and integrate over $\Omega$ obtaining 
\begin{align} 
c_v \intO{ \vr |D_t \vt|^2 } & + \frac{\kappa}{2} \frac{\D }{\dt} \intO{ |\Grad \vt |^2 } \br &= 
c_v \intO{ \vr D_t \vt \ \vu \cdot \Grad \vt } - \intO{ \vr \vt \ \Div \vu \ D_t \vt} + \intO{ \vr \vt \  \Div \vu\  \vu \cdot \Grad \vt } \br 
&+ \frac{\D }{\dt} \intO{ \vt \  \mathbb{S}(\Ds \vu) : \Grad \vu } \br	&- 
\mu \intO{ \vt \left( \Grad \vu + \Grad^t \vu - \frac{2}{3} \Div \vu \mathbb{I} \right): \left( \Grad \partial_t \vu + \Grad^t \partial_t \vu - \frac{2}{3} \Div \partial_t \vu \mathbb{I} \right) } \br
&- 2 \eta \intO{ \vt \  \Div \vu \  \Div \partial_t \vu }. 	
\label{g12}
	\end{align}
Indeed the term involving the boundary integral is handled as 
\[
- \kappa \intO{ \Del \vt \ \partial_t \vt } = - \kappa \int_{\partial \Omega} \partial_t \vtB \Grad \vt \cdot \vc{n} \ \D S_x + \frac{\kappa}{2}  \frac{\D }{\dt} \intO{ |\Grad \vt |^2 }, 
\]
where
\[
\int_{\partial \Omega} \partial_t \vtB \Grad \vt \cdot \vc{n} \ \D S_x = 0
\]
as the boundary temperature is independent of $t$.

Similarly to Fang, Zi, Zhang \cite[Section 3, Lemma 3.4]{FaZiZh}, we have to show that the intergrals 
\[
\intO{ \vt \  \Grad \vu : \Grad \partial_t \vu },\ 
\intO{ \vt \ \Grad \vu : \Grad^t \partial_t \vu },\ \mbox{and} \ \intO{ \vt \ \Div \vu \ \Div \partial_t \vu }
\]
can be rewritten in the form compatible with \eqref{g4}, meaning with the time derivatives replaced by material 
derivatives. Fortunately, this step can be  carried out in the present setting using only the boundary condition $\vu \cdot \vc{n}|_{\partial \Omega} = 0$. Indeed  we get 
\[
\intO{ \vt \ \Grad \vu : \Grad \partial_t \vu } = 
\intO{ \vt \ \Grad \vu : \Grad (D_t \vu) } - \intO{ \vt \ \Grad \vu : \Grad (\vu \cdot \Grad \vu)},
\]
where
\begin{align*}
	&\intO{ \vt \  \Grad \vu : \Grad (\vu \cdot \Grad \vu)} \br
	&= \intO{\vt \ \Grad \vu : (\Grad \vu \cdot \Grad \vu)}  + \frac{1}{2} \intO{ \vt \ \vu \cdot \Grad | \Grad \vu |^2 } \br 
	&= 
	\intO{\vt \ \Grad \vu : (\Grad \vu \cdot \Grad \vu)} - \frac{1}{2} \intO{ |\Grad \vu|^2 \ \Grad \vt \cdot \vu } 
	- \frac{1}{2} \intO{ |\Grad \vu|^2 \  \vt \Div \vu }.
\end{align*}

Similarly, 
\[
\intO{ \vt \  \Grad \vu : \Grad^t \partial_t \vu } = 
\intO{ \vt \ \Grad \vu : \Grad^t (D_t \vu) } - \intO{ \vt \ \Grad \vu : \Grad^t (\vu \cdot \Grad \vu)},
\]
where 
\begin{align}
&\intO{ \vt \ \Grad \vu : \Grad^t (\vu \cdot \Grad \vu)} \br  
&=  \intO{\vt \ \Grad \vu : (\Grad^t \vu \cdot \Grad^t \vu)} + 
\frac{1}{2} \intO{ \vt \ \vu \cdot \Grad ( \Grad \vu : \Grad^t \vu ) } \br 
&= \intO{\vt \ \Grad \vu : (\Grad^t \vu \cdot \Grad^t \vu)}- 
\frac{1}{2} \intO{ ( \Grad \vu : \Grad^t \vu ) \ \Grad \vt \cdot \vu }- \frac{1}{2} \intO{ ( \Grad \vu : \Grad^t \vu ) \  \vt  \Div \vu }.
\nonumber
\end{align}
Finally, 
\[
\intO{ \vt \ \Div \vu \ \Div \partial_t \vu } = 
\intO{ \vt \ \Div \vu \ \Div  D_t \vu } - \intO{ \vt \ \Div \vu \ \Div (\vu \cdot \Grad \vu) },
\]
where 
\begin{align} 
&\intO{ \vt \ \Div \vu \ \Div (\vu \cdot \Grad \vu) }  \br 
&= \intO{ \vt \ \Div \vu \ (\Grad \vu : \Grad^t \vu)} + \frac{1}{2} \intO{ \vt \vu \cdot \Grad |\Div \vu |^2 } \br 
&= \intO{ \vt \ \Div \vu \ (\Grad \vu : \Grad^t \vu)} - \frac{1}{2} \intO{ |\Div \vu |^2 \ \Grad \vt \cdot \vu }-\frac{1}{2} \intO{ |\Div \vu |^2 \  \vt \Div \vu    }.
	\nonumber
\end{align}

We conclude, using \eqref{g2}, \eqref{g4}, and \eqref{g12}, 
\begin{align} 
\intO{ |\Grad \vt |^2 (\tau, \cdot) } &+ \int_0^\tau \intO{ \vr |D_t \vt |^2 } \dt \br &\leq c(T, \mathcal{D}_0, \Ov{\vr}, \Ov{\vt}, \Ov{u} ) \left(1 	
+ \int_0^\tau \intO{ |\Grad \vu |^4 } \dt \right). \label{g13}
	\end{align}

Next, by virtue of the decomposition $\vu = \vc{v} + \vc{w}$ and the bound \eqref{g9}, 
\begin{equation} \label{g14}
\intO{ |\Grad \vu |^4 } \aleq \intO{ |\Grad \vc{v} |^4 } + \intO{ |\Grad \vc{w} |^4 }  
\leq c(T, \mathcal{D}_0, \Ov{\vr}, \Ov{\vt}, \Ov{u} ) \left( 1 + \intO{ |\Grad \vc{w} |^4 } \right), 
\end{equation}
and, similarly, 
\begin{equation} \label{g15}
\| \vc{w} \|_{L^\infty (\Omega; R^3) } \leq \| \vc{u} \|_{L^\infty (\Omega; R^3) } + \| \vc{v} \|_{L^\infty (\Omega; R^3) } \leq c(T, \mathcal{D}_0, \Ov{\vr}, \Ov{\vt}, \Ov{u} ) .
\end{equation}
Recalling the Gagliardo--Nirenberg interpolation inequality in the form 
\begin{equation} \label{g16}
\| \Grad U \|_{L^4(\Omega; R^3)}^2 \leq \|  U \|_{L^\infty(\Omega)} \|  \Del  U \|_{L^2(\Omega)} \ \mbox{whenever}\ U|_{\partial \Omega} = 0, 
\end{equation}
we may use \eqref{g14}, \eqref{g15} to rewrite \eqref{g13} in the form 
\begin{align} 
	\intO{ |\Grad \vt |^2 (\tau, \cdot) } &+ \int_0^\tau \intO{ \vr |D_t \vt |^2 } \dt \br &\leq c(T, \mathcal{D}_0, \Ov{\vr}, \Ov{\vt}, \Ov{u} ) \left(1 + \int_0^\tau \intO{ |\Grad \vt |^2 } \dt	
	+ \int_0^\tau \| \vc{w} \|_{W^{2,2}(\Omega; R^3)}^2 \dt \right). \label{g17}
\end{align}

Finally, we use the elliptic estimates \eqref{g11} to conclude
\begin{align} 
	&\intO{ |\Grad \vt |^2 (\tau, \cdot) } + \int_0^\tau \intO{ \vr |D_t \vt |^2 } \dt \br &\leq c(T, \mathcal{D}_0, \Ov{\vr}, \Ov{\vt}, \Ov{u} ) \left(1 + \int_0^\tau \intO{ \left( |\Grad \vt |^2 
	+ |\Grad \vu |^2 \right) } \dt	
	+ \int_0^\tau \| \sqrt{\vr} \partial_t \vu \|_{L^2(\Omega; R^3)}^2 \dt \right). \label{g18}
\end{align}
Summing up \eqref{g2}, \eqref{g4}, and \eqref{g18} we may apply Gronwall's lemma to obtain the following bounds:
\begin{align} 
\sup_{t \in [0,T)} \| \vu (t, \cdot) \|_{W^{1,2}(\Omega; R^3)} &\leq c(T, \mathcal{D}_0, \Ov{\vr}, \Ov{\vt}, \Ov{u} ), \label{g19}\\
\sup_{t \in [0,T)} \|\sqrt{\vr}  D_t \vu (t, \cdot) \|_{L^{2}(\Omega; R^3)} &\leq c(T, \mathcal{D}_0, \Ov{\vr}, \Ov{\vt}, \Ov{u} ), \label{g20} \\
\sup_{t \in [0,T)} \| \vt (t, \cdot) \|_{W^{1,2}(\Omega)} &\leq c(T, \mathcal{D}_0, \Ov{\vr}, \Ov{\vt}, \Ov{u} ), \label{g21}\\ 
\int_0^T \intO{ |\Grad D_t \vu |^2 } \dt &\leq c(T, \mathcal{D}_0, \Ov{\vr}, \Ov{\vt}, \Ov{u} ), \label{g22}	\\
\int_0^T \intO{ \vr |D_t \vt |^2 } \dt &\leq c(T, \mathcal{D}_0, \Ov{\vr}, \Ov{\vt}, \Ov{u} ). \label{g23}
	\end{align}
Moreover, it follows from \eqref{g9}, \eqref{g16}, \eqref{g20} 
\begin{equation} \label{g24}
	\sup_{t \in [0,T) } \| \Grad \vu (t, \cdot) \|_{L^4(\Omega; R^{3\times 3})} \leq c(T, \mathcal{D}_0, \Ov{\vr}, \Ov{\vt}, \Ov{u} ).
\end{equation}
In addition, \eqref{g23}, \eqref{g24} and the standard parabolic estimates applied to the internal energy balance \eqref{i3} yield 
\begin{equation} \label{g25}
	\int_0^T \| \vt \|^2_{W^{2,2}(\Omega)} \dt \leq c(T, \mathcal{D}_0, \Ov{\vr}, \Ov{\vt}, \Ov{u} ).
	\end{equation}

\section{Second energy bound}
\label{s}

It follows from \eqref{g11}, \eqref{g20} that
\begin{equation} \label{s1}
	\sup_{t \in [0,T) } \| \vc{w}(t, \cdot) \|_{W^{2,2}(\Omega; R^3)} \leq c(T, \mathcal{D}_0, \Ov{\vr}, \Ov{\vt}, \Ov{u} ); 
\end{equation}
whence, by virtue of \eqref{g9} and Sobolev embedding $W^{1,2}(\Omega) \hookrightarrow L^6(\Omega)$, 
\begin{equation} \label{s2}
	\sup_{t \in [0,T)} \| \Grad \vu (t, \cdot) \|^2_{L^6(\Omega; R^{3\times 3})}  \leq c(T, \mathcal{D}_0, \Ov{\vr}, \Ov{\vt}, \Ov{u} ).
\end{equation}
Moreover, as a consequence of \eqref{g22}, $D_t \vu$ is bounded in $L^2(L^6)$, which, combined with \eqref{s2}, gives rise to 
\begin{equation} \label{s3}
\int_0^T \left\| \partial_t \vu \right\|^2_{L^6(\Omega; R^3)} \dt\leq c(T, \mathcal{D}_0, \Ov{\vr}, \Ov{\vt}, \Ov{u} ). 
\end{equation}
Finally, going back to \eqref{g7} we conclude 
\begin{equation} \label{s4}
\int_0^T \| \vc{w} \|_{W^{2,6}(\Omega; R^3)}^2 \dt \leq c(T, \mathcal{D}_0, \Ov{\vr}, \Ov{\vt}, \Ov{u}),
\end{equation}
and 
\begin{equation} \label{s5}
\int_0^T \| \vc{u} \|_{W^{1,q}(\Omega; R^3)}^2 \dt \leq c(T, \mathcal{D}_0, \Ov{\vr}, \Ov{\vt}, \Ov{u},q ) \ \mbox{for any} \ 1 \leq q < \infty.	
\end{equation}

\section{Estimates of the derivatives of the density}
\label{d}

Using \eqref{s4}, \eqref{s5}, we may proceed as in \cite[Section 5]{SuWaZh} to deduce the bounds
\begin{equation} \label{d1}
	{\rm sup}_{t \in [0,T)} \left( \| \partial_t \vr (t, \cdot) \|_{L^6(\Omega)} + \| \vr (t,\cdot) \|_{W^{1,6} (\Omega)} \right) \leq c(T, \mathcal{D}_0, \Ov{\vr}, \Ov{\vt}, \Ov{u} ).
\end{equation}
Revisiting the momentum equation \eqref{i2} we use \eqref{d1} together with the other bounds established above to obtain 
\begin{equation} \label{d2}
	\int_0^T \| \vu \|^2_{W^{2,6}(\Omega; R^3)} \dt \leq c(T, \mathcal{D}_0, \Ov{\vr}, \Ov{\vt}, \Ov{u} ).
\end{equation}

\subsection{Positivity of the density and temperature}

It follows from \eqref{d2} that $\Div \vu$ is bounded in $L^1(0,T; L^\infty(\Omega))$. Thus the equation of continuity \eqref{i1} yields a positive lower bound on the density 
\begin{equation} \label{d3}
 \inf_{(t,x) \in [0,T)\times\Ov{\Omega}}\ \vr(t,x) \geq \underline{\vr} > 0, 
\end{equation}
where the lower bound  depends on the data as well as on the length $T$ of the time interval.

Similarly,
rewriting the internal energy balance equation \eqref{i3} in the form 
\begin{equation} \label{d4}
c_v \left( \partial_t \vt + \vu \cdot \Grad \vt \right) - \frac{\kappa}{\vr} \Del \vt =  \frac{1}{\vr} \mathbb{S} : \Ds \vu - \vt \Div \vu  
\end{equation}
we may apply the standard parabolic maximum/minimum principle to deduce 
\begin{equation} \label{d5}
	 \inf_{(t,x) \in [0,T)\times\Ov{\Omega}} \ \vt(t,x) \geq \underline{\vt} > 0. 
\end{equation}

\section{Parabolic regularity for the heat equation}
\label{p}

We rewrite the parabolic equation \eqref{d4} in terms of $\Theta = \vt - \vtB$. 
Recalling $\Del \vtB = 0$ we get 
\begin{equation} \label{p2}
	c_v \left( \partial_t \Theta + \vu \cdot \Grad \vt \right) - \frac{\kappa}{\vr} \Del \Theta =  \frac{1}{\vr} \mathbb{S} : \Ds \vu - \vt \Div \vu  
\end{equation}
with the \emph{homogeneous} Dirichlet boundary conditions
\begin{equation} \label{p3}
\Theta|_{\partial \Omega} = 0.	
	\end{equation} 

Now, we can apply all arguments of \cite[Sections 4.6, 4.7]{FeSu2015_N1} to $\Theta$ obtaining the bounds 
\begin{align} 
	\| \vt \|_{C^\alpha([0,T] \times \Ov{\Omega})} &\leq c(T, \mathcal{D}_0, \Ov{\vr}, \Ov{\vt}, \Ov{u} )
	\ \mbox{for some}\ \alpha > 0, \label{p4} \\
\| \vt \|_{L^p(0,T; W^{2,3}(\Omega))} + 
\| \partial_t \vt \|_{L^p(0,T; L^{3}(\Omega))}  &\leq  c(T, \mathcal{D}_0, \Ov{\vr}, \Ov{\vt}, \Ov{u} )
\ \mbox{for all}\ 1 \leq p < \infty, \label{p5}
\end{align}
together with 
\begin{equation} \label{p6}
\| \vu \|_{L^p(0,T; W^{2,6}(\Omega;R^3))} + 
\| \partial_t \vu \|_{L^p(0,T; L^{6}(\Omega;R^3))}  \leq  c(T, \mathcal{D}_0, \Ov{\vr}, \Ov{\vt}, \Ov{u} )
\ \mbox{for any}\ 1 \leq p < \infty.
\end{equation}

\section{Final estimates}
\label{f}

The bounds \eqref{p6} imply, in particular, 
\begin{equation} \label{f1}
\sup_{(t,x) \in [0,T)\times \Ov{\Omega}}|\Grad \vu (t,x) | \leq c(T, \mathcal{D}_0, \Ov{\vr}, \Ov{\vt}, \Ov{u} ).	
\end{equation}	
Thus the desired higher order estimates can be obtained exactly as in \cite[Section 4.6]{FeNoSun1}. Indeed the arguments of \cite[Section 4.6]{FeNoSun1} 
are based on differentiating the equation \eqref{p2} with respect to time which gives rise to a parabolic problem for $\partial_t \vt$ with the \emph{homogeneous} Dirichlet boundary 
conditions $\partial_t \vt|_{\partial \Omega} = 0$. Indeed we get 
\begin{equation*}
\begin{aligned}
	{ c_v \partial_{tt}^2 \vt  + c_v \vu \cdot \Grad \partial_t \vt - \frac{\kappa}{\vr} \Delta_x \partial_t\vt =} & { - c_{v} \partial_t \vu \cdot \Grad \vt  - \frac{1}{\vr^2}\partial_t \vr \left( \kappa \Delta_x \vt + \mathbb{S} (\mathbb{D}_x \vu): \mathbb{D}_x \vu \right)} \\
	&{ + \frac{2}{\vr}\ \mathbb{S} (\mathbb{D}_x\vu): \mathbb{D}_x \partial_t\vu-\partial_t \vt \ \Div \vu - \vt \ \Div \partial_t \vu.}
\end{aligned}
\end{equation*} 
The estimates obtained in the previous sections imply that the right--hand side of the above equation 
is bounded in $L^2(0,T; L^2(\Omega))$.
Thus multiplying the equation on $\Del \partial_t \vt$ and performing the standard by parts integration, we get 
the desired estimates as in \cite[Section 4.6]{FeNoSun1}.

The remaining estimates are obtained exactly as in \cite[Section 4.6]{FeNoSun1} :
\begin{align} \label{f2}
	\sup_{t \in [0,T) } \| \vt (t, \cdot) \|_{W^{3,2}(\Omega)} + 	\sup_{t \in [0,T) } \| \partial_t \vt (t, \cdot) \|_{W^{1,2}(\Omega)} &\leq  c(T, \mathcal{D}_0, \Ov{\vr}, \Ov{\vt}, \Ov{u} ), \\
	\int_0^T \left( \| \partial_t \vt \|^2_{W^{2,2}(\Omega)} + \| \vt \|^2_{W^{4,2}(\Omega)} \right) \dt &\leq c(T, \mathcal{D}_0, \Ov{\vr}, \Ov{\vt}, \Ov{u} ), \label{f3} \\
	\sup_{t \in [0,T) } \| \vu (t, \cdot) \|_{W^{3,2}(\Omega; R^3)} + 	\sup_{t \in [0,T) } \| \partial_t \vu (t, \cdot) \|_{W^{1,2}(\Omega; R^3)} &\leq c(T, \mathcal{D}_0, \Ov{\vr}, \Ov{\vt}, \Ov{u} ) , \label{f4} \\
		\int_0^T \left( \| \partial_t \vu \|^2_{W^{2,2}(\Omega; R^3)} + \| \vu \|^2_{W^{4,2}(\Omega; R^3)} \right) \dt &\leq c(T, \mathcal{D}_0, \Ov{\vr}, \Ov{\vt}, \Ov{u} ), \label{f5} 
	\end{align}
and 
\begin{equation} \label{f6}
	\sup_{t \in [0,T) } \| \vr (t, \cdot) \|_{W^{3,2}(\Omega)} \leq c(T, \mathcal{D}_0, \Ov{\vr}, \Ov{\vt}, \Ov{u} ).
	\end{equation} 	

We have completed the proof of Proposition \ref{PT}. 

\def\cprime{$'$} \def\ocirc#1{\ifmmode\setbox0=\hbox{$#1$}\dimen0=\ht0
	\advance\dimen0 by1pt\rlap{\hbox to\wd0{\hss\raise\dimen0
			\hbox{\hskip.2em$\scriptscriptstyle\circ$}\hss}}#1\else {\accent"17 #1}\fi}



\end{document}